\documentclass[11pt]{article}
\input amssymb.sty
\input amssym.def
\input amssym
\input epsf

\newtheorem{theorem}{Theorem}[section]

\newtheorem{claim}[theorem]{Claim}
\newtheorem{theorem-definition}[theorem]{Theorem-Definition}
\newtheorem{theorem-construction}[theorem]{Theorem-Construction}
\newtheorem{lemma-definition}[theorem]{Lemma--Definition}
\newtheorem{proposition-definition}[theorem]{Proposition--Definition}
\newtheorem{lemma}[theorem]{Lemma}
\newtheorem{proposition}[theorem]{Proposition}

\newtheorem{definition}[theorem]{Definition}

\begin{document}
\newcommand{\Z}{{\mathbb Z}}
\newcommand{\R}{{\mathbb R}}
\newcommand{\Q}{{\mathbb Q}}
\newcommand{\C}{{\mathbb C}}
\newcommand{\lms}{\longmapsto}
\newcommand{\lra}{\longrightarrow}
\newcommand{\hra}{\hookrightarrow}
\newcommand{\ra}{\rightarrow}
\newcommand{\sgn}{\rm sgn}
\begin{titlepage}
\title{Cluster ensembles, quantization  and  the dilogarithm II: 
The intertwiner}  
\author{V.V. Fock, A. B. Goncharov}
\end{titlepage}
\date{\it To Yuri Ivanovich Manin for his 70th birthday}
\maketitle

\tableofcontents

 \vskip 6mm \noindent

\section{Introduction}

{\bf Abstract.} This paper is the second part of our preprint 
``Cluster ensembles, quantization  and  the dilogarithm'' \cite{FG2}. 
Its main result is a construction, by means of the quantum dilogarithm, 
 of certain intertwiner operators, which play the crucial role in the 
quantization of the cluster ${\cal X}$-varieties and construction of the corresponding canonical representation. 

\vskip 3mm

Cluster ensemble ({\it loc. cit.}) is a pair $({\cal A}, {\cal X})$ of schemes over $\Z$, called 
cluster ${\cal A}$- and ${\cal X}$-varieties, related by a map $p: {\cal A} \lra {\cal X}$. 
The ring of regular functions on the cluster ${\cal A}$-variety is the upper cluster algebra 
\cite{BFZ}. 

Cluster ${\cal A}$- and ${\cal X}$-varieties are glued from families of coordinatized split algebraic tori 
by means of certain subtraction free rational transformations. In particular 
it makes sense to consider the spaces of their positive real points, denoted by ${\cal A}^+$ and ${\cal X}^+$.

The cluster ${\cal X}$-variety has a Poisson structure, given in any cluster coordinate system 
$\{X_i\}$ 
by 
$$
\{X_i, X_j\}= \widehat \varepsilon_{ij}X_iX_j, \quad \widehat \varepsilon_{ij}\in \Z. 
$$
The Poisson tensor $\widehat \varepsilon_{ij}$ depends on the choice of coordinate system. 
There is a canonical non-commutative deformation ${\cal X}_q$ of the cluster ${\cal X}$-variety 
in the direction of this Poisson structure (\cite{FG2}).

\vskip 3mm 
The gluing procedure underlying the definition of a cluster variety can be 
understood as a functor from a certain groupoid, called the {\it cluster modular groupid} ${\cal G}$, 
to a category of commutative algebras. 

 In Section 3 we suggest a $\ast$-quantization 
of the cluster ${\cal X}$-variety, understood as a functor from the gropoid ${\cal G}$ 
to the category of non-commutative topological $\ast$-algebras. 
More precisely, the coordinate systems on cluster varieties 
are parametrised by the objects ${\bf i}$ of the groupoid ${\cal G}$,  called {\it seeds}. 
To each seed ${\bf i}$ we assign two coordinatized 
tori, ${\cal A}_{\bf i}$ and ${\cal X}_{\bf i}$. The algebra of 
the smooth functions on the latter admits a canonical $\hbar$-deformation, given by a topological 
Heisenberg
 $\ast$-algebra ${\cal H}^{\hbar}_{\bf i}$. 
So to define a functor we need to relate these algebras for different seeds. 
We write the formulas relating the generators 
of the algebras ${\cal H}^{\hbar}_{\bf i}$, but do not specify 
the category of topological 
$\ast$-algebras. As a result, 
the $\ast$-quantization 
of the cluster ${\cal X}$-variety serves only 
as a motivation, and  we 
state in Sections 3 Claims instead of Theorems when 
those unspecified topological algebras enter the formulations. 
Hopefully there will be a precise version of Section 3. 
However the rest of the paper does not depend on that, while motivations 
given in Section 3 
clarify what we do next. 

\vskip 3mm 
In Section 4 we proceed to a construction of the canonical unitary projective representation 
of the modular gropoid. It is realized in the 
Hilbert space $L^2({\cal A}^+)$ assigned to the set of positive real points of the 
cluster ${\cal A}$-variety.

For each seed ${\bf i}$ there is a Hilbert space $L^2({\cal A}^+_{\bf i})$ (which is  canonically identified 
with $L^2({\cal A}^+)$). In Section 4.1 the Heisenberg $\ast$-algebra ${\cal H}^{\hbar}_{\bf i}$ 
is represented in $L^2({\cal A}^+_{\bf i})$. 
In fact a bigger algebra -- the chiral double ${\cal H}^{\hbar}_{\bf i}\otimes {\cal H}^{\hbar}_{\bf i^o}$ of the Heisenberg $\ast$-algebra -- acts on the same Hilbert space.  

The morphisms in the groupoid ${\cal G}$ are defined as compositions of certain elementary ones, called 
{\it mutations} and {\it symmetries}. 
Given a mutation ${\bf i}\to {\bf i'}$ we construct a unitary operator 
$$
{\bf K}_{{\bf i}\to {\bf i'}}: L^2({\cal A}^+_{\bf i}) \lra L^2({\cal A}^+_{\bf i'}).  
$$
It intertwines the actions of 
the Heisenberg $\ast$-algebras   related to the seeds ${\bf i}$ and ${\bf i'}$\footnote{The precise meaning we put into 
 ``intertwining'' is clear from the computation carried out in the proof: we deal with the genrators 
of the Heisenbrg algebra only, and thus are not concerned with the nature of the topological completion, which 
was left  unspecified in Section 3. For a different approach see \cite{FG3}.}. 
(Similar itertwiners for symmetries are rather tautological). 
This construction is the main result of the paper. 
The operator ${\bf K}_{{\bf i}\to {\bf i'}}$ is characterised by its intertwining property uniquely up to a constant. 
 
Certain compositions of mutations and symmetries are identity morphisms 
in the  groupoid ${\cal G}$. So to get a representation of the modular groupoid 
we have to show that the corresponding compositions of the intertiners ${\bf K}_{{\bf i}\to {\bf i'}}$ 
are multiples of the identity operators. This is a rather difficult problem. It is solved 
in \cite{FG3}, where we give   another construction of the intertwiner 
and introduce the geometric object reflecting its 
properties, the cluster double. Alltogether the intertwiners 
give rise to a unitary projective reprersentation 
of the cluster modular groupoid.

\vskip 3mm
The paper is organised as follows: the essential for us 
 properties of the quantum logarithm and dilogarithm are collected, without proofs in Section 5. The proofs and more 
of the properties of these functions can be found in Section 4 of \cite{FG3}. 
In Section 2.1 we recall,  for the convinience of the reader, 
basic definitions/facts  about cluster ensembles. 
 Claim \ref{11.12.03.323} delivers 
a quantization of the space of real positive points 
of the cluster ${\cal X}$-variety. 
The main results of this paper is Theorem \ref{11.11.03.76} providing an explicit formula for the 
intertwiner 
${\bf K}_{{\bf i}\to {\bf i'}}$. 

\vskip 3mm
{\bf Acknowledgments}. 
V.F. was supported by the grants CRDF 2622; 2660. 
A.G. was 
supported by the  NSF grants  DMS-0099390 and DMS-0400449. 

\section{Cluster ensembles} 

\subsection{Basic definitions} A {\em seed} ${\mathbf i}$ is a triple $(I, \varepsilon, d)$, 
where $I$ is a finite set, 
$\varepsilon$ is a matrix $\varepsilon_{ij}$, $i,j \in I$, with 
$\varepsilon_{ij} \in {\mathbb Z}$, and $d = \{d_i\}$, $i \in I$, 
are positive integers, such that the matrix 
$\widehat{\varepsilon}_{ij}:=\varepsilon_{ij}d^{-1}_j$ is skew-symmetric.

For a seed ${\mathbf i}$ 
we assign a torus ${\mathcal X}_{\mathbf i} = ({\mathbb G}_m)^I$ with  the 
coordinates $\{X_i| i\in I\}$ on the factors. 
and a Poisson structure given by 
\begin{equation} \label{f1}
\{X_i,X_j\}=\widehat{\varepsilon}_{ij}X_iX_j.
\end{equation} 

Let ${\mathbf i}=(I, \varepsilon, d)$ and ${\mathbf i}'=(I', \varepsilon', d')$  be two seeds, and $k\in
I$. A {\em mutation in the direction $k\in I$} is an isomorphism $\mu_k: I\rightarrow I'$ such that $d'_{\mu_k(i)}=d_i$, and 
\begin{equation} \label{f2}
\varepsilon'_{\mu_k(i)\mu_k(j)}=
\left\{ \begin{array}{lll} -\varepsilon_{ij} & \mbox{ if }  i=k \mbox{ or } j=k,\\
\varepsilon_{ij}& \mbox{ if } \varepsilon_{ik}\varepsilon_{kj} \leq  0,\\
\varepsilon_{ij} + |\varepsilon_{ik}|\varepsilon_{kj} & \mbox{ if } \varepsilon_{ik}\varepsilon_{kj}> 0.
\end{array}\right.
\end{equation} 
A {\em symmetry} of a seed ${\mathbf i}=(I, \varepsilon, d)$ is an  automorphism $\sigma$ of the set $I$ preserving the matrix $\varepsilon$ and the numbers $d_i$. 
Symmetries and mutations induce rational maps between the corresponding  
seed ${\cal X}$-tori, denoted by the same symbols $\mu_k$ and $\sigma$ and given by the formulae 
$\sigma^*X_{\sigma(i)}=X_i$ 
and
\begin{equation} \label{f3}
\mu_{k(i)}^*X_{\mu_k(i)} = \left\{\begin{array}{lll} X_k^{-1}& \mbox{ if } & i=k, \\
    X_i(1+X_k^{-\sgn(\varepsilon_{ik})})^{-\varepsilon_{ik}} & \mbox{ if } &  i\neq k. \\
\end{array} \right.
\end{equation} 

A {\em seed cluster transformation} 
is a composition of symmetries and mutations. 
Two seeds are {\em equivalent} if related by a cluster transformation. The equivalence class of a seed ${\mathbf i}$ is denoted by $|{\mathbf i}|$.
A seed cluster transformation induces a rational map
between the two seed
${\cal X}$-tori, called a cluster transformation map. 

A {\em cluster ${\cal X}$-variety ${\cal X}_{|{\mathbf i}|}$} is a scheme over $\Z$ obtained by gluing the
 seed ${\cal X}$-tori  for the seeds equivalent to a given 
seed ${\mathbf i}$  
via the cluster transformation maps, and then taking 
the affine closure. 
   Every seed 
  provides a cluster ${\cal X}$-variety with a rational  coordinate system. 
Its coordinates are called {\em cluster coordinates}. 
Cluster transformation maps preserve the Poisson structure. 
Thus a cluster ${\cal X}$-variety
 has a canonical 
Poisson structure. 

\vskip 3mm 
 {\em The cluster ${\cal A}$-varieties}. Given a seed ${\bf i}$, define 
a {\em seed ${\cal A}$-torus} ${\mathcal A}_{\mathbf i}: = ({\mathbb G}_m)^I$ 
with the standard  
standard coordinates $\{A_i| i\in I\}$ on the factors. 
Symmetries and mutations give rise to  birational 
maps between the 
seed ${\cal A}$-tori, given by 
$\sigma^*A_{\sigma(i)}=A_i$ 
and
\begin{equation} \label{f3cc}
\mu_{k(i)}^*A_{\mu_k(i)} = \left\{\begin{array}{lll} A_i
 & \mbox{ if } & i\not =k, \\
    A_k^{-1}\left(\prod_{i|\varepsilon_{ki}>0}A_i^{\varepsilon_{ki}} + 
\prod_{i|\varepsilon_{ki}<0}A_i^{-\varepsilon_{ki}}\right) & \mbox{ if } &  i= k. \\
\end{array} \right.
\end{equation} 
The cluster ${\cal A}$-variety ${\cal A}_{|{\mathbf i}|}$ 
is a scheme over $\Z$ obtained by gluing 
 all seed ${\cal A}$-tori for the seeds equivalent to a given 
seed ${\mathbf i}$ 
  using the above birational isomorphisms, and taking 
the affine closure.

\vskip 3mm
There is a map 
$
p: {\cal A} \lra {\cal X}, 
$
given in every cluster coordinate system by $p^*X_k = \prod_{I\in I}A_i^{\varepsilon_{ki}}$. 
\vskip 3mm

Cluster ${\cal A}$- and ${\cal X}$-varieties have canonical positive atlases, so it makes sense to consider the sets of their 
real positive points, denoted ${\cal A}^+$ and ${\cal X}^+$.

\vskip 3mm
{\it The cluster modular groupoid}. 
Seed cluster transformations inducing the same map of the seed ${\cal A}$-tori 
are called {\it trivial seed cluster transformations}. 
The {\it cluster modular groupoid} ${\cal G}_{|{\bf i}|}$ 
is a groupoid whose objects are 
seeds equivalent to a given seed ${\bf i}$, and  
$Hom({\bf i}, {\bf i'})$ is the set of all 
seed cluster transformations from ${\bf i}$ to  ${\bf i'}$ 
modulo the trivial ones. 
Given  a seed ${\bf i}$, 
the {\it cluster mapping class group $\Gamma_{\bf i}$} is the automorphism group of the object ${\bf i}$ of ${\cal G}_{|{\bf i}|}$. 
The group $\Gamma_{\bf i}$
 acts by automorphisms 
of the cluster ${\cal A}$-variety. 

\vskip 3mm
{\it The quantum space ${\cal X}_q $}. It is a canonical non-commutative 
$q$-deformation of the cluster ${\cal X}$-variety 
defined in Section 3 of \cite{FG2}. 
We start from the {\it seed quantum torus algebra} ${\rm T}^q_{\bf i}$, defined as 
an associative $\ast$-algebra with generators 
$X^{\pm 1}_i$, $i\in I$ and $q^{\pm 1}$ and relations
$$
q^{-\widehat \varepsilon_{ij}}X_iX_j = q^{-\widehat \varepsilon_{ji}}X_jX_i,\qquad 
\ast X_i = X_i, \quad \ast q = q^{-1}.
$$
Let ${\rm QTor}^*$ be a category whose objects are 
quantum torus algebras and morphisms are 
$\ast$-homomorphisms of their fraction fields. 
The quantum space ${\cal X}_q $ is understood 
as a contravariant functor 
$$
\eta^q: 
\mbox{The  modular groupoid ${\cal G}_{|{\bf i}|}$} \lra {\rm QTor}^*.
$$
It assigns to 
a seed ${\bf i}$ the quantum torus $\ast$-algebra ${\rm T}^q_{{\bf i}}$, and 
to a mutation   ${\bf i} \lra {\bf i'}$ a map of the fraction fields 
${\rm Frac}({\rm T}^q_{{\bf i'}}) \lra {\rm Frac}({\rm T}^q_{{\bf i}})$, given by 
a $q$-deformation of formulas (\ref{f3})\footnote{See also a more transparent and 
less computational
 definition given later on in Section 3 of \cite{FG3}.}: Set $q_k:= q^{1/h_k}$, 
$$
(\mu_k^q)^*(X_i):=  \left\{\begin{array}{lll} X_k^{-1}
 & \mbox{ if } & i=k, \\
    X_i\left(\prod_{a =1}^{|\varepsilon_{ik}|} (1+q_k^{2a-1}X^{-\sgn (\varepsilon_{ik})}_k)\right)^{-\varepsilon_{ik}} & \mbox{ if } &  i\not = k. \\
\end{array} \right. 
$$
One uses Theorem 7.2 in {\it loc. cit.}  to prove that 
$\eta^q$ sends trivial seed cluster transformations to the identity maps.

\vskip 3mm
The {\it chiral dual} to a 
seed ${\bf i} = (I, \varepsilon, d)$ is  a  
seed 
${\bf i}^o:= (I, -\varepsilon, d)$. Mutations commute with the chiral duality on seeds. Therefore a cluster 
${\cal X}$-variety ${\cal X}$ (respectively ${\cal A}$-variety ${\cal A}$), 
gives rise to the 
chiral dual cluster 
${\cal X}$-variety (respectively ${\cal A}$-variety) denoted by 
 ${\cal X}^o$ (respectively ${\cal A}^0$). 
They are related, see Lemmas \ref{6.9.03.11} - \ref{6.9.03.11s} . 

The {\it Langlands dual} to a seed 
${\bf i} = (I, \varepsilon, d)$ is the seed 
 ${\bf i}^{\vee} = (I, \varepsilon^{\vee}, d^{\vee})$, where 
$d_i^{\vee}: = d_i^{-1}$ 
and 
\begin{equation} \label{6.4.03.2}
\varepsilon^{\vee}_{ij} = -\varepsilon^{\vee}_{ji} := \widehat d^{-1}_i\varepsilon_{ij}\widehat d_j, \qquad 
\widehat d_i: = d_i^{-1}.
\end{equation}
The Langlands duality on seeds commutes with mutations. Therefore 
it gives rise to the Langlands dual cluster 
${\cal A}$-, and ${\cal X}$-varieties, denoted 
${\cal A}^{\vee}$ and ${\cal X}^{\vee}$.

\vskip 3mm
Below we skip the subscript 
$|{\bf i}|$ encoding  
the corresponding 
 cluster ensemble whenever possible.

\vskip 3mm
\subsection{Connections between quantum ${\cal X}$-varieties}
 There are three ways to alter the space ${\cal X}_{|{\bf i}|, q}$: 

(i) change $q$ to $q^{-1}$, 

(ii) change ${\bf i}$ to ${\bf i}^o$, 

(iii) change the quantum space ${\cal X}_{|{\bf i}|, q}$ to the opposite quantum space 
${\cal X}^{\rm opp}_{|{\bf i}|, q}$.

\noindent
(In (iii) we change 
every quantum torus from which we glue the space to the opposite one).

The following lemma tells that the resulting three quantum  spaces are canonically isomorphic:
\begin{lemma} \label{6.9.03.11} 
a) There is a canonical isomorphism of quantum spaces
$$
\alpha_{\cal X}^q: {\cal X}_{|{\bf i}|, q} \lra {\cal X}^{\rm opp}_{|{\bf i}|, q^{-1}}, \qquad (\alpha_{\cal X}^q)^*
: X_i \lms X_i. 
$$
(Given in on the generators of any cluster coordinate system by $X_i \lms X_i$).

b) There is a canonical isomorphism of quantum spaces
$$
i^q_{\cal X}: {\cal X}_{|{\bf i}|, q} \lra {\cal X}_{|{\bf i}^o|, {q^{-1}}}, \qquad (i^q_{\cal X})^*: X_i^o\lms X_i^{-1}.
$$
(Given in any cluster coordinate system by $X^o_i \lms X_i^{-1}$, where 
$X^o_i$ are the generators of ${\cal X}_{{\bf i}^o, {q^{-1}}}$).

c) There is a canonical isomorphism of quantum spaces
$$
\beta_{\cal X}^q:= \alpha_{\cal X}^q\circ i^q_{\cal X}: 
{\cal X}_{|{\bf i^o}|, q} \lra {\cal X}^{\rm opp}_{|{\bf i}|, q}, \qquad X_i \lms {X_i^o}^{-1}. 
$$
\end{lemma}

{\bf Proof}. Apparently each of the three maps is an isomorphism of the corresponding seed quantum tori algebras. 
For example, in the case b) we have 
$$
i_{\cal X}^q \Bigl((q^{-1})^{-\widehat {\varepsilon}^o_{ij}} X^o_i
X^o_j \Bigr) = q^{-\widehat\varepsilon_{ij}}  X_i^{-1} X_j^{-1}.
$$
So we need to check that they commute with the mutations. 

a) 
Let us assume first that $\varepsilon_{ik} = a<0$. The claim results from the fact
 that the following two compositions are equal (observe that $\alpha^*$ is an antiautomorphism):
$$
X_i \stackrel{\mu_k^*}{\lms}  X_i \prod_{b=1}^a(1+q^{2b-1}_kX_k) \stackrel{\alpha^*}{\lms} \prod_{b=1}^a(1+q^{2b-1}_kX_k) X_i; 
$$
$$
X_i \stackrel{\alpha^*}{\lms} X_i \stackrel{(\mu^o_k)^*}{\lms}  X_i \prod_{b=1}^a(1+q^{-(2b-1)}_kX_k). 
$$
The computation in the case $\varepsilon_{ik} >0$ is similar. 

b) 
To check that 
$i_{\cal X}^q\circ \mu_{k}^* = 
\mu_{k}^*\circ i_{\cal X}^q$ we  calculate 
each of the  maps on the generator $X_i$. 
Let us assume $\varepsilon_{ik} = -a < 0$. Then 
$\varepsilon^o_{ik} = a$, and one has 
$$
i_{\cal X}^q\circ \mu_{k}^*({X^o_i}') = 
X_i^{-1} \prod_{b=1}^a(1+q^{2b-1}X^{-1}_k),
$$
$$
\mu_{k}^*\circ i_{\cal X}^q({X^o_i}') = 
\Bigl(X_i \prod_{b=1}^a(1+q^{-(2b-1)}X^{-1}_k)^{-1}\Bigr)^{-1} = 
\prod_{b=1}^a(1+q^{-(2b-1)}X^{-1}_k) X_i^{-1} = 
X_i^{-1}\prod_{b=1}^a(1+q^{2b-1}X^{-1}_k). 
$$
The case $\varepsilon_{ik} > 0$ is similar.  One can deduce it to the case $\varepsilon_{ik} < 0$ 
since $\mu_k \circ \mu_k = {\rm Id}$, and $\varepsilon_{ik}' = - \varepsilon_{ik}$. The part b) is proved. 

c) Follows from a) and b). The lemma is proved.

\vskip 3mm
\begin{lemma} \label{6.9.03.11s} 
The cluster 
ensembles  related to the seeds ${\bf i}$ and ${\bf i^o}$ are  
canonically isomorphic as pairs of varieties.   The isomorphism is provided by the 
following maps:
$$
{\rm Id}: {\cal A}_{|{\bf i}|} \lra {\cal A}_{{\bf i^o}}; 
\qquad \quad i_{\cal X}: {\cal X}_{|{\bf i}|} \lra {\cal X}_{|{\bf i^o}|}.
$$
\end{lemma}

{\bf Proof}. In a given cluster coordinate system our maps are obviously 
isomorphisms. The 
compatibility with ${\cal X}$-cluster transformations 
 is the part b) of  Lemma \ref{6.9.03.11}; for the ${\cal A}$-cluster transformations  we have 
$$
A^o_k {A^o}'_k~ \stackrel{}{=}
~\prod_{\varepsilon^o_{ki}>0} 
(A^o_i)^{\varepsilon^o_{ki}}
+ \prod_{\varepsilon^o_{ki}<0}  
(A^o_i)^{-\varepsilon^o_{ki}} ~= ~
\prod_{\varepsilon_{ki}<0} A_i^{- \varepsilon_{ki}}
+ \prod_{ \varepsilon_{ki}>0} 
A_i^{ \varepsilon_{ki}}~ \stackrel{}{=} ~ A_k A'_k 
$$
Compatibility with the projection $p$ is clear. The lemma is proved. 
\vskip 3mm

\section{Motivation: $\ast$-quantization of cluster ${\cal X}$-varieties}

\subsection{$\ast$-quantization of the space ${\cal X}^+$ via the quantum logarithm}

Let $\{X_i\}$ be coordinates on the cluster ${\cal X}$-variety
 corresponding to a seed ${\bf i}$. 
Since by definition the functions $X_i$ are strictly positive at the 
points of ${\cal X}^+$, we can introduce  
the {\it logarithmic coordinates} $
x_i := \log X_i 
$ on ${\cal X}^+$. 
For every seed ${\bf i}$ they provide an isomorphism 
$$
\beta_{\bf i}: 
 {\cal X}_{\bf i}^+ \stackrel{\sim}{\lra} \R^I; \qquad t \lms \{x_i(t)\}.
$$
For a mutation $\mu_k: {\bf i} \to {\bf i}'$ there is a gluing map 
\begin{equation} \label{5.24.03.1as}
\beta_{{\bf i} \to {\bf i}'}: {\cal X}^+_{\bf i} \lra {\cal X}^+_{\bf i'}, 
\qquad \beta_{{\bf i} \to {\bf i}'}(x'_i) =  
\left\{ \begin{array}{ll} 
  x_i  - \varepsilon_{ik} \log (1+e^{- {\rm sgn}(\varepsilon_{ik}) x_k})&  i \not = k, 
\\
-x_k & i  = k\end{array}\right.
\end{equation}
\vskip 3mm

To prepare the soil for quantization, let us look at this from 
a different point  of view.  Let ${\rm Com}^*$ be the category of 
commutative topological $\ast$-algebras over $\C$.  
Recall the cluster modular groupoid ${\cal G}_{|{\bf i}|}$. 
There is a 
contravariant functor
$$
\beta: {\cal G}_{|{\bf i}|} \lra {\rm Com}^*. 
$$
Namely, we assign to a seed ${\bf i}$ a commutative topological $\ast$-algebra 
$S({\cal X}^+_{\bf i})$ of smooth complex valued 
functions in ${\cal X}^+_{\bf i}$ with 
$\ast f:= \overline f$, and to a mutation ${\bf i} \to {\bf i'}$ 
a homomorphism $\beta^*_{{\bf i} \to {\bf i}'}: S({\cal X}^+_{\bf i'}) 
\lra S({\cal X}^+_{\bf i})$.

\vskip 3mm

Let ${\cal C}$ be a category whose morphisms are 
$\C$-vector spaces. {\it Projectivisation} $P{\cal C}$ 
of the category  ${\cal C}$ as a 
 new category with the same objects as ${\cal C}$, and morphisms given by  
${\rm Hom}_{P{\cal C}}(C_1, C_2):= {\rm Hom}_{{\cal C}}(C_1, C_2)/U(1)$, where $U(1)$ is the multiplicative group of complex numbers with absolute value $1$. 
A {\it projective functor} 
$F: {\cal G} \to {\cal C}$ is a functor from ${\cal G}$ to $P{\cal C}$. 

Let ${\cal C}^{\ast}$ be the category of  topological $\ast$-algebras. 
Two functors  $F_1, F_2: {\cal C} \lra {\cal C}^{\ast}$ 
{\it essentially coincide} if there exists a third functor $F$ and natural transformations $F_1 \to F, F_2 \to F$ providing 
for every object $C$  dense inclusions  
$F_1(C) \hra F_3(C)$, $F_1(C) \hra F_3(C)$.

\begin{definition} \label{6.3.03.10} 
A {\it quantization of  the space ${\cal X}^+_{|{\bf i}|}$} 
is a family of contravariant projective functors 
$$
\kappa_{|{\bf i}|}^\hbar: {\cal G}_{|{\bf i}|} \lra {\cal C}^{\ast}
$$
depending smoothly on a real parameter $\hbar$, 
related to the original Poisson 
manifold ${\cal X}^+_{|{\bf i}|}$ as follows:

i) The limit  $\kappa_{|{\bf i}|}:= \lim_{\hbar \to 0}\kappa_{|{\bf i}|}^\hbar$ exists 
and essentially coincides with the functor $\beta$ defining 
${\cal X}^+_{|{\bf i}|}$.

ii) The Poisson bracket given by 
$\lim_{\hbar \to 0}[f_1,f_2]/\hbar$ is defined  and 
coincides with the one on ${\cal X}^+_{|{\bf i}|}$. 
\end{definition}

Let us define a quantization functor $\kappa^\hbar = \kappa^\hbar_{|{\bf i}|}$. 
We assign to every seed  ${\bf i}$  the 
Heisenberg $\ast$--algebra ${\cal H}^\hbar_{{\bf i}}$. 
It is a topological 
$\ast$--algebra over $\C$ generated 
by the elements $x_i$ such that 
$$
[x_j, x_k] = 2 \pi i \hbar \widehat \varepsilon_{jk}; \quad 
x_j^* = x_j;\qquad q = e^{\pi i \hbar}.
$$

Further, let us  assign to mutation $\mu_k: {\bf i} \to {\bf i'}$ 
a homomorphism of topological 
$\ast$-algebras 
\begin{equation} \label{6.1.03.2}
\kappa^\hbar(\mu_k): {\cal H}^\hbar_{{\bf i'}} \lra {\cal H}^\hbar_{{\bf i}}.
\end{equation}
We  employ the quantum logarithm $\phi^\hbar(z)$, see 
(\ref{phi}). 
 Denote by $x_i'$ the generators of ${\cal H}^\hbar_{{\bf i'}}$. 
 Set
\begin{equation} \label{6.1.03.1}
\hbar_k:= \widehat d_k \hbar; \qquad\quad \kappa^\hbar(\mu_k): x_i' \lms  
\left\{ \begin{array}{ll}x_i -  
 \varepsilon_{ik}\phi^{\hbar_k}(-{\rm sgn}(\varepsilon_{ik}) x_k) &  
\mbox{ if $k\not =i$},\\
-x_i & \mbox{ if $k=i$}. 
\end{array}\right.
 \end{equation}

\begin{claim} \label{3.13.03.1} a)  Formulas (\ref{6.1.03.1}) 
provide a morphism  of $\ast$--algebras  
$\kappa_{|{\bf i}|}^\hbar(\mu_k): {\cal H}^\hbar_{{\bf i'}} 
\lra {\cal H}^\hbar_{{\bf i}}$. 

b) 
The  collection of 
 $\ast$-algebras $\{{\cal H}_{\bf i}\}$ 
and morphisms 
$\{\kappa_{|{\bf i}|}^\hbar({\mu_k})\}$ provide a quantization functor 
$$
\kappa_{|{\bf i}|}^\hbar: 
{\cal G}_{|{\bf i}|} \lra P{\cal C}^*.
$$

c) Let $\hbar^{\vee}:= 1/\hbar$. Then there are  isomorphisms 
\begin{equation} \label{6.3.03.1}
{\cal H}^\hbar_{{\bf i}} \lra {\cal H}^{\hbar^{\vee}}_{{\bf i}^{\vee}}, 
\qquad x_i \lms x_i^{\vee}:= \frac{x_i}{\widehat d_i \hbar}; \qquad 
 \end{equation}
They give rise to  a natural transformation of functors 
$\kappa^\hbar_{|{\bf i}|} \lra \kappa^{\hbar^{\vee}}_{|{\bf i}^{\vee}|}$. 
\end{claim}

{\bf Justification}. a) 
Property {\bf A3} of the function $\phi^{\hbar}(x)$, see
 Section 4, guarantees that the morphism 
$\kappa^\hbar(\mu_k)$ preserves the real structure. 
It follows from Property {\bf A1} that when $\hbar \to 0$ the limit 
of the quantum formula 
(\ref{6.1.03.1}) exists and coincides with the mutation formula  
(\ref{5.24.03.1as}) 
for the logarithmic coordinates $x_i$. 

 b) To check that we have a functor one needs to check first that 
mutation formulas are compatible with the 
transformations $\kappa^\hbar(\mu_k)$. This is a straitforward calculation
using  Property {\bf A5}. Then
one has to check that the defining relations for the groupoid 
${\cal G}_{|{\bf i}|}$ are mapped to zero. 
Here we need the results of Sections 3.2-3.3 of \cite{FG2} 
and the following well known lemma:

\begin{lemma} \label{6.4.03.10}
Suppose that $A,B$ are selfadjoint operators,
 $[A,B] =-\lambda$ is a scalar, and $f(z)$ is a continuous  
function with primitive $F(z)$. Then 
$$
e^{A+f(B)} = e^A {\rm exp}\left\{\frac{1}{\lambda}\int^{B+\lambda}_{B}f(z)dz\right\}: = 
e^A {\rm exp}\left(\frac{F(B+\lambda) - F(B)}{\lambda}\right).
$$
\end{lemma}

c)  Thanks to 
formula (\ref{6.4.03.2}) the map $x_i \lms x^{\vee}_i$ is 
an $\ast$--algebra homomorphism: 
$$
[x_i^{\vee}, x_j^{\vee}]~~ = ~~\frac{[x_i, x_j]}{\widehat d_i\widehat d_j \hbar^2}~~ = ~~
2\pi i \hbar^{\vee}\widehat \varepsilon_{ij}/\widehat d_i\widehat d_j~~ \stackrel{(\ref{6.4.03.2})}{=} ~~
2\pi i \hbar^{\vee}\widehat \varepsilon^{\vee}_{ij}.
$$
To verify that it commutes with 
mutation homomorphisms we use Properties {\bf A2} and {\bf A4} 
of the function $\phi^{\hbar}(x)$, observing that 
$$
\frac{x_i+ |\varepsilon_{ik}| \phi^{\hbar_k}(x_k)}{\widehat d_i \hbar}~~  = ~~
 \frac{x_i}{\widehat d_i \hbar} + 
\frac{|\widehat d_i^{-1}\varepsilon_{ik}\widehat d_k|\phi^{\widehat d_k \hbar}(x_k)}
{\widehat d_k \hbar}~~  \stackrel{(\ref{6.4.03.2}) + A4}{=}  ~~ x^{\vee}_i + |\varepsilon^{\vee}_{ik}| 
\phi^{\hbar^{\vee}_k}(x^{\vee}_k).
$$

\vskip 3mm

\subsection{Modular double of a cluster ${\cal X}$-variety and 
$\ast$-quantization of  the space ${\cal X}^+$} 
Set
$$
q:= e^{\pi i\hbar}, \quad q^{\vee}:= e^{\pi i/\hbar}, \quad \hbar \in \R.
$$ 
 
\begin{definition} \label{6.2.03.3} 
The modular double  ${\cal X}_{|{\bf i}|, q}\times 
{\cal X}_{|{\bf i}^{\vee}|, q^{\vee}}$ of a quantum cluster ${\cal X}$-variety 
${\cal X}_{|{\bf i}|, q}$ 
is  a contravariant functor  
$$
\eta_{|{\bf i}|}^q\otimes \eta_{|{\bf i}^{\vee}|}^{q^{\vee}}: 
{\cal G}_{|{\bf i}|} \lra 
{\rm QTor}^*.
$$
\end{definition}
\vskip 3mm
So we assign to a seed ${\bf i}$ 
a quantum torus algebra
 ${\rm T}_{{\bf i}}^q \otimes {\rm T}_{{\bf i}^{\vee}}^{q^{\vee}}$, and  
to a mutation $\mu_k: {\bf i}\to {\bf i'}$ a positive  
$\ast$-homomorphism of the fraction fields of the quantum torus algebras 
$$
\eta_{|{\bf i}|}^q(\mu_k)\otimes \eta_{|{\bf i}^{\vee}|}^{q^{\vee}}(\mu_k): 
{\Bbb T}_{{\bf i'}}^q \otimes 
{\Bbb T}_{{\bf i'}^{\vee}}^{q^{\vee}} 
\lra {\Bbb T}_{{\bf i}}^q \otimes 
{\Bbb T}_{{\bf i}^{\vee}}^{q^{\vee}}, \qquad {\Bbb T}:= {\rm Frac}({\rm T}).
$$

We want to relate the modular double  ${\cal X}_{|{\bf i}|, q}\times 
{\cal X}_{|{\bf i}^{\vee}|, {q^{\vee}}}$ with the quantization of the  space ${\cal X}^+_{|{\bf i}|}$. 
We are going to define a natural transformation of functors
$ \eta_{|{\bf i}|}^q\otimes \eta_{{\cal E}^{\vee}}^{q^{\vee}} 
\lra \kappa^\hbar_{|{\bf i}|}$. 

\vskip 3mm
We use the following easy fact. Assume that 
$[y_i,y_j]$ is a scalar. Then we have 
\begin{equation} \label{6.2.03.1}
e^{y_i}e^{y_j} = e^{[y_i,y_j]}e^{y_j}e^{y_i}.
 \end{equation}
Let ${\bf i}$ be a seed. 
Denote by $X_i$ the generators of ${\Bbb T}_{\bf i}^q$, and 
by $X^{\vee}_i$ the generators of ${\Bbb T}_{{\bf i}^{\vee}}^{q^{\vee}}$. 
It is easy to check using (\ref{6.2.03.1}) 
that there are the following homomorphisms: 
$$
l_{{\bf i}}: {\Bbb T}_{{\bf i}}^q \lra {\cal H}_{{\bf i}}^h, \quad X_i \lms e^{x_i}, \quad \mbox{and} \quad 
 l^{\vee}_{\bf i}: 
{\Bbb T}_{{\bf i}^{\vee}}^{q^{\vee}} \lra {\cal H}_{{\bf i}}^{\hbar}, 
\quad X^{\vee}_i \lms 
e^{x^{\vee}_i}.
$$

They evidently commute 
with the $\ast$-structures. Their images commute. Indeed, since 
$\widehat \varepsilon_{ij} \in \Z$ one has $e^{[x_i, x_j/\hbar]} = 
e^{2\pi i \widehat \varepsilon_{ij}} =1$. So 
  $e^{x_i}$ commutes with $e^{x_j/\hbar}$. Therefore 
they  give rise to a  homomorphism of the tensor product:
$$
{L}_{\bf i}:= l_{\bf i}\otimes  l^{\vee}_{\bf i}: {\Bbb T}_{{\bf i}}^q 
\otimes {\Bbb T}_{{\bf i}^{\vee}}^{q^{\vee}}\lra {\cal H}_{{\bf i}}^\hbar.
$$

\begin{proposition} \label{erer}
For any mutation $\mu_k: {\bf i} \to {\bf i'}$ 
the following diagram, where the left vertical arrow is the map 
$ \eta^q(\mu_k)\otimes \eta^{q^{\vee}}(\mu_k)$, and the right one is $\kappa^\hbar(\mu_k)$, 
is   commutative: 
 $$
\begin{array}{ccc}
{\Bbb T}_{\bf i}^q \otimes {\Bbb T}_{\bf i}^{q^{\vee}}& 
\stackrel{L_{\bf i}}{\lra}& {\cal H}_{\bf i}^\hbar\\
&&\\
\eta^{q, q^{\vee}}\uparrow &
&\uparrow \kappa^\hbar\\
&&\\
{\Bbb T}_{\bf i'}^q \otimes 
{\Bbb T}_{\bf i'}^{q^{\vee}}&\stackrel{L_{\bf i'}}{\lra}&
{\cal H}_{\bf i'}^\hbar
\end{array}
$$
\end{proposition}

{\bf Proof}. We need Lemma \ref{6.4.03.10}. 
The case $i = k$ is trivial, so we assume that $i \not = k$. 
Let  $\varepsilon_{ik} = -a \leq 0$. 
Then applying the lemma we get
$$
\kappa_{|{\bf i}|}^\hbar(\mu_k) L_{{\bf i'}} (X'_i \otimes 1)~ = ~\kappa^\hbar(\mu_k) e^{x^{\vee}_i} ~ = ~ e^{x_i + a\phi^{\hbar_k}(x_k)}~ =  ~\frac{e^{x_i}}{2\pi i a \hbar_k}
{\rm exp}\left(\int_{x_k}^{x_k+ 2 \pi i \hbar_k a}a\phi^{\hbar_k}(z)dz\right)~ = 
$$
$$
~\frac{e^{x_i} }{2\pi i \hbar_k}
{\rm exp}\left(\int_{- \infty}^{x_k}\Bigl(\phi^{\hbar_k}(z+2 \pi i \hbar_k a ) - 
\phi^{\hbar_k}(z)\Bigr)dz\right)~ \stackrel{A5}{=}
~\frac{e^{x_i} }{2\pi i \hbar_k}
{\rm exp}\left(\int_{- \infty}^{x_k}\sum^{a}_{b=1}
\frac{2\pi i \hbar_k}{e^{-z-i\pi (2b-1)\hbar_k}+1}dz\right)~ = 
$$
$$ 
e^{x_i} \prod_{b=1}^{a}(1+q_k^{2a-1}e^{x_k})~ = ~L_{{\bf i}}
\left(X_i \prod_{b=1}^{a}(1+q_k^{2a-1}X_k)\right)~
=  ~L_{{\bf i}}(\eta_{|{\bf i}|}^q (\mu_k)
\otimes \eta_{|{\bf i}^{\vee}|}^{q^{\vee}}(\mu_k))(X'_i \otimes 1).
$$
The calculation in the case $\varepsilon_{ik} = a \geq 0$ is similar. 
 The proposition is proved.
\vskip 3mm

\begin{claim} \label{11.12.03.323} The collection 
of homomorphisms $\{L_{\bf i}\}$ 
provides a morphism of functors 
\begin{equation} \label{6.2.03.4}
{\Bbb L}^\hbar: \eta_{|{\bf i}|}^q\otimes 
\eta_{|{\bf i^{\vee}|}}^{q^{\vee}} \lra \kappa_{|{\bf i}|}^\hbar.
\end{equation}
\end{claim}

{\bf Justification}. Is given by Proposition \ref{erer}. 

\vskip 3mm
{\bf Representations of the quantized ${\cal X}^+_{|{\bf i}|}$-space.}  
The following definition serves as a motivation of the 
construction of intertwiners presented below. 

\begin{definition} \label{11.12.03.1}
A projective $\ast$-representation of the 
quantized ${\cal X}^+_{|{\bf i}|}$-space is the following data:

i) A projective functor 
$$
{\cal L}_{|{\bf i}|}: {\cal G}_{|{\bf i}|} \lra \mbox{the category of Hilbert spaces}.
$$
It includes for each object ${\bf i}$ of ${\cal G}_{|{\bf i}|}$ 
a Hilbert spaces ${\Bbb L}_{\bf i}$, 
and for every mutation $\mu_k: {\bf i} \to {\bf i'}$ a unitary operator,
defined up to a scalar of absolute value $1$:
$$
{\bf K}_{{\bf i},{\bf i'}}: {\Bbb L}_{\bf i} \lra {\Bbb L}_{\bf i'}. 
$$

ii) A $\ast$-representation $\rho_{\bf i}$ 
of the Heisenberg algebra ${\cal H}^\hbar_{\bf i}$ 
in the Hilbert space ${\Bbb L}_{\bf i}$.

iii) The operators ${\bf K}_{{\bf i}, {\bf i'}}$ 
intertwine the representations ${\rho}_{\bf i}$ 
and ${\rho}_{{\bf i'}}$:
$$
{\rho}_{\bf i}(s) = {\bf K}_{{\bf i} ,{\bf i'}}^{-1} {\rho}_{{\bf i'}}\Bigl(\kappa^\hbar
(\mu_k)(s)\Bigr) {\bf K}_{{\bf i},{\bf i}'}, \qquad s \in {\cal H}^\hbar_{{\bf i}}.
$$
\end{definition}
The morphisms of the representations of the quantum ${\cal X}^+_{|{\bf i}|}$-space 
are defined in an obvious way. 

\vskip 3mm
{\it Representations of the 
mapping class group $\Gamma_{|{\bf i}|}$}. 
Restricting the functor $\rho_{|{\bf i}|}$ to the group of automorphisms of 
an object of the groupoid ${\cal G}_{|{\bf i}|}$ we get a projective unitary representation of 
$\Gamma_{|{\bf i}|}$. 

\vskip 3mm
The Heisenberg algebra ${\cal H}^\hbar_{{\bf i}}$ has a 
family of irreducible $\ast$-representation by operators in a Hilbert space. 
These representations are characterized by the central character $\chi$.

\vskip 3mm
The collection of the Hilbert spaces 
$\{{\Bbb L}_{{\bf i}}\}$ and representations $\{\rho_{{\bf i}}\}$ 
is by no means canonical: it depends, for example, on the 
choice of polarization of the Heisenberg algebra. Once choosen, 
it determines the intertwiners ${\bf K}_{{\bf i},{\bf i}'}$. 
Below we introduce a {\it canonical representation} 
of the  chiral double of the quantized space ${\cal X}_{|{\bf i}|}^+$, defined 
by using the Hilbert spaces 
$L^2({\cal A}^+_{\bf i})$. 
\vskip 3mm

\section{The intertwiner}

\subsection{A bimodule structure on functions 
on  the ${\cal A}$--space} 
Let $X$ be an algebra. Recall that  $M$ is a bimodule 
over $X$ if $X$ acts on $M$ from the left as well as from the right, 
and these two actions  commute. So $M$ is 
an $X \otimes X^{\rm opp}$-module, where 
 $X^{\rm opp}$ is the algebra with the product $x\ast y:= y x$.

\vskip 3mm

Let us choose a seed ${\bf i}$. Recall the algebra  $\Q[{\cal A}_{\bf i}]$  
of regular functions on the seed torus ${\cal A}_{\bf i}$. 
We assume that $q\in \C^*$. 
For each $i \in I$ let us define commuting algebra homomorphisms
$$
t^{\pm}_i: \Q[{\cal A}_{\bf i}] \lra \Q[{\cal A}_{\bf i}]; \quad t^{\pm}_i: \left\{ \begin{array}{ll} 
A_i \lms q^{^{\pm}\widehat d_i}A_i& \\
A_j \lms A_j & j \not = i.
\end{array}\right.
$$
Since $\varepsilon_{ii} =0$,  $A_i$ 
does not appear in the monomial $p^*X_i$, and so 
the operator of multiplication by $p^*X_i$ commutes with 
$t_i^{\pm}$. 
Let us define an ${T}^q_{\bf i}$-bimodule structure on $\Q[{\cal A}_{\bf i}]$. 
The left and right actions of the generator $X_i$ on $f \in \Q[{\cal A}_{\bf i}]$ are given by 
\begin{equation}\label{5.29.03.2}
X_i \circ f:= p^*X_i \cdot t^-_i(f), \qquad f\circ X_i 
:= p^*X_i \cdot t^+_i(f),
\end{equation}
\begin{lemma} \label{5.29.03.1}
The operators (\ref{5.29.03.2}) provide $\Q[{\cal A}_{\bf i}]$ with a structure  of a bimodule 
over the algebra ${T}^q_{\bf i}$. 
\end{lemma}

{\bf Proof}. Observe that one has 
$$
t_j^+(p^*X_i)\cdot p^*X_j 
= t_i^-(p^*X_j)\cdot p^*X_i = 
q^{\widehat \varepsilon_{ij}}p^*X_ip^*X_j.
$$
Indeed, the first term equals to 
$ q^{\varepsilon_{ij}d_j} p^*X_i\cdot p^*X_j $, and the second 
is  $ q^{-\varepsilon_{ji}d_i} p^*X_i\cdot p^*X_j $. 
One has 
$$
q^{-\widehat \varepsilon_{ij}}X_iX_j \circ f = 
q^{-\widehat \varepsilon_{ij}}p^*X_i \cdot t^-_i(p^*X_j) \cdot t^-_it^-_j(f) = 
p^*X_i \cdot p^*X_j \cdot t^-_it^-_j(f). 
$$
$$
f\circ q^{-\widehat \varepsilon_{ij}}X_iX_j  = 
q^{-\widehat \varepsilon_{ij}}p^*X_j \cdot t^+_j(p^*X_i) \cdot t^+_jt^+_i(f) = 
p^*X_i \cdot p^*X_j \cdot t^+_it^+_j(f). 
$$
Since the right hand sides are evidently symmetric in $i,j$ 
we have the desired relations, and hence the left and right 
actions of the quantum algebra torus.  
Further, the two
actions commute:
$$
X_i \circ (f \circ X_j) = X_i \circ \Bigl(t_j^+(f) p^* X_j\Bigr) = p^*X_i t_i^-(
p^* X_j) t_i^-t_j^+(f).
$$
$$
(X_i \circ f) \circ X_j = (p^* X_i t_i^-(f))\circ  X_j = p^*X_j t_j^+(p^* X_i) 
t_j^+t_i^-(f).
$$
The lemma is proved. \vskip 3mm

\vskip 3mm
{\it The logarithmic version of the bimodule structure}.  
Since the coordinate functions  $A_i$ are positive on the 
space ${\cal A}^+$, 
  one can introduce new coordinates $a_j:= \log A_j$. 
They provide an isomorphism 
$
\alpha_{\bf i}: {\cal A}_{\bf i}^+ \stackrel{\sim}{\lra} \R^{I}. 
$
Set $da:= da_1 \wedge \cdots  \wedge da_{|I|}$.  
There is a Hilbert space 
$L^2({\cal A}_{\bf i}^+)$ with a scalar product $$
 (f, g):= 
\int_{{\cal A}_{\bf i}^+}f(a)\overline {g(a)}da.
$$
Apparently the form $da$ changes the sign under a  mutation ${\bf i} \to {\bf i'}$. 
So the Hilbert spaces $L^2({\cal A}_{\bf i}^+)$ for different seeds 
${\bf i}$ are naturally identified. 
Consider the following operators in $L^2({\cal A}_{\bf i}^+)$:
$$
\widehat x_j^-:=  -\pi i \hbar \widehat d_j \frac{\partial}{\partial a_j} +
\sum_k\varepsilon_{jk}a_k, 
\qquad 
\widehat x_j^+:=  \pi i  \hbar \widehat d_j \frac{\partial}{\partial a_j} + 
\sum_k\varepsilon_{jk}a_k. 
$$
\begin{lemma} \label{5.29.03.5}
The 
operators $\{\widehat x_j^{\pm}\}$  provide the 
Hilbert space $L^2({\cal A}_{\bf i}^+)$ 
with a structure of a bimodule over the $\ast$--algebra 
${\cal H}^\hbar_{\bf i}$.  
\end{lemma}

{\bf Proof}. These  operators are selfadjoint and one has 
$$
[\widehat x_j^-, \widehat x_k^-] = 2 \pi i \hbar \widehat \varepsilon_{jk}; \qquad 
[\widehat x_j^+, \widehat x_k^+] = -2 \pi i \hbar \widehat \varepsilon_{jk}, 
\qquad 
[\widehat x_j^-, \widehat x_k^+]=0 \quad \mbox{for any $j,k \in I$}.
$$
The lemma is proved.
\vskip 3mm 
{\bf Remark}. There is an automorphism $x_j \lms - x_j$ 
of the Heisenberg algebra 
${\cal H}_{\bf i} $. Similarly there is 
an automorphism $X_j \lms  X^{-1}_j$ of the quantum torus algebra $T^q_{\bf i}$.

\vskip 3mm

\subsection{The intertwiner via the quantum dilogarithm }
 Let $\mu_k: {\bf i} \to {\bf i}'$ be a mutation. 
By Lemma \ref{5.29.03.5} for each seed ${\bf i}$ 
the Hilbert 
space $L^2({\cal A}^+_{\bf i})$ has a natural 
${\cal H}^\hbar_{\bf i}$-bimodule structure.  According to Lemma \ref{6.9.03.11}, 
this is the same as the ${\cal H}_{{\bf i}}^h\otimes {\cal H}_{{\bf i^o}}^\hbar$-module structure. 
Our goal is to  define an operator
\begin{equation} \label{11.13.03.2}
{\bf K}_{{\bf i} \to {\bf i}'}: L^2({\cal A}^+_{\bf i}) \lra L^2({\cal A}^+_{\bf i'}) 
\end{equation}
intertwining the ${\cal H}_{{\bf i}}^\hbar\otimes {\cal H}_{{\bf i^o}}^\hbar$- and 
${\cal H}_{{\bf i'}}^\hbar\otimes {\cal H}_{{{\bf i'}^o}}^\hbar$-module structures. 
By this we mean only that the operator ${\bf K}_{{\bf i} \to {\bf i}'}$ intertwines the action of the generators 
of ${\cal H}_{{\bf i}}^\hbar\otimes {\cal H}_{{\bf i^o}}^\hbar$ withthe action of their images 
under the gluing map map $\kappa_{|{\bf i}|}^{\hbar}({{\bf i} \to {\bf i}'})$. 
\vskip 3mm

{\it The function $G$}. 
Let us introduce our key function
$$
G(a_1,\ldots,a_n)=
$$
\begin{equation} \label{11.13.03.1}
\int \Phi^{\hbar_k}(\widehat{d}_kc
-\sum_j\varepsilon_{kj}a_j)^{-1}\Phi^{\hbar_k}(-\widehat{d}_kc
-\sum_j\varepsilon_{kj}a_j)\exp\left(c
\frac{\sum\limits_{j|\varepsilon_{kj}<0}\varepsilon_{kj}a_j+
a_k}{\pi i \hbar}\right)dc.
\end{equation}
Substituting the explicit integral expression for the function
$\Phi^{\hbar_k}(z)$ one gets
$$
G(a_1,\ldots,a_n)=\int\exp
\left(\int_{\Omega}\frac{\exp(it\sum_j
\varepsilon_{kj}a_j)\sin(t \widehat{d}_kc)}{2it~{\rm sh}(\pi t)
{\rm sh}(\pi \hbar_kt)}  dt + c\frac{\sum\limits_{j|
\varepsilon_{kj}<0}\varepsilon_{kj}a_j+a_k}{\pi i \hbar } \right)dc.
$$

We denote by $(a_1, ..., a_n)$ the logarithmic coordinates 
corresponding to $v$, and by $(a_1, ..., a_k', ..., a_n)$ 
the ones corresponding to ${\bf i}'$. Recall that only the 
coordinate $a_k$ changes under the mutation $\mu_k$. 
Let us define the operator ${\bf K}_{{\bf i} \to {\bf i}'}$ by 
\begin{equation}\label{convolution}
({\bf K}_{{\bf i} \to {\bf i}'}f)(a_1,\ldots, a'_k, \ldots, a_n) := \int G(a_1,\ldots, 
a'_k+a_k, \ldots, a_n)f(a_1,\ldots, a_k,\ldots, a_n)da_k,
\end{equation}
where $a'_k+a_k$ and $a_k$ are on the $k$-th places.

\begin{theorem} \label{11.11.03.76}
The operators 
${\bf K}_{{\bf i} \to {\bf i}'}$ 
 intertwine the ${\cal H}^\hbar_{{\bf i}}\otimes {\cal H}^\hbar_{{\bf i^o}}$-module 
structures on $L^2({\cal A}^+_{\bf i})$ provided by Lemma \ref{5.29.03.5}. 
\end{theorem}

\vskip 3mm
{\bf Remark}. We prove in \cite{FG3} that the collection of Hilbert space 
$L^2({\cal A}^+_{\bf i})$ and operators ${\bf K}_{{\bf i} \to {\bf i}'}$ 
provide a unitary projective representation 
of the groupoid ${\cal G}_{|{\bf i}|}$. 
This implies that 
the operators ${\bf K}_{{\bf i} \to {\bf i}'}$ give rise 
to a unitary projective representation of the cluster modular group 
$\Gamma_{|\bf i|}$ in $L^2({\cal A}^+_{\bf i})$. 
\vskip 3mm

{\bf Proof}. 
We present a computation which allows to find the   
function ${G}$ as a unique up to a scalar 
function such that the corresponding 
integral transformation   intertwines the ${\cal H}^\hbar_{\bf i}$- and 
${\cal H}^\hbar_{{\bf i'}}$-bimodule structures on $L^2({\cal A}^+_{\bf i})$ and 
$L^2({\cal A}^+_{\bf i'})$.
Recall that $\varepsilon^o_{ij} = -\varepsilon_{ij}$, so we may write
$\varepsilon^{\pm}_{ij} := \pm  \varepsilon_{ij}$ and denote by $x_i^{\pm}$
the $x_i$-coordinates for the seeds ${\bf i}$ and ${\bf i^o}$. 

So we have to find ${G}$ such that the 
integral transformation (\ref{convolution}) induces a map of operators:
\begin{equation} \label{12.9.03.10}
{\widehat x^{'\pm}}_i \lms \left\{ \begin{array}{ll}{\widehat x}^{ \pm}_i -
\varepsilon_{ik}^{\pm}\phi^{\hbar_k}(-{\rm sgn}(\varepsilon^{\pm}_{ik})\widehat 
x_k^{\pm})  &  
\mbox{ if $i \not = k$}\\
-\widehat x^{\pm}_k & \mbox{ if $i=k$}. 
\end{array}\right.
\end{equation}
This means that we should have (changing $\widehat x^-_i$ to $-\widehat x^-_i$ 
for convenience) for $i \not = k$: 
$$
\pi i h\widehat{d_i} \frac{\partial}{\partial a_i} \pm \sum_{j\neq k}
\varepsilon'_{ij}a_j \pm \varepsilon'_{ik}a'_k
\mapsto
$$
\begin{equation}\label{homomorphism}
\mapsto
\pi i \hbar \widehat{d}_i\frac{\partial}{\partial a_i} \pm
\sum_{j\neq k} \varepsilon_{ij}a_j \pm \varepsilon_{ik}\left(a_k -
\phi^{\hbar_k}\left(-\sgn(\pm \varepsilon_{
ik})(\pi i \hbar \widehat{d}_k\frac{\partial}{\partial a_k} \pm
\sum_j\varepsilon_{kj}a_j)\right)\right),
\end{equation}
and
\begin{equation}\label{homomorphism_k}
\pi i \hbar \widehat{d}_k \frac{\partial}{\partial a'_k}  \pm
\sum_j\varepsilon'_{kj}a_j \mapsto 
-\Bigl(\pi i \hbar \widehat{d}_k \frac{\partial}{\partial a_k}  \pm
\sum_j\varepsilon_{kj}a_j \Bigr).
\end{equation}
Here we use the following conventions. The signs $\pm$ in our formulas 
always use either $+$ everywhere, or $-$ everywhere, so $\mp:= -\pm$.
Thus we have one set of the equations corresponding to
the upper signs and another one  to the lower signs.

Observe that $\varepsilon_{kk}=0$. The relation (\ref{homomorphism_k}) 
is satisfied by (\ref{convolution})
if and only if 
$$- \varepsilon'_{kj} = \varepsilon_{kj}, \quad \mbox{and} \quad  
\frac{\partial}{\partial a'_k} \lms -\frac{\partial}{\partial a_k}.
$$ 
Since these two conditions are evidently valid, we have 
the relation (\ref{homomorphism_k}). 

Substituting (\ref{homomorphism}) into (\ref{convolution}) one gets
the identities:
$$
\int(\pi i \hbar  \widehat{d_i}\frac{\partial G}{\partial a_i}f + \pi i \hbar \widehat{d_i}
G\frac{\partial f}{\partial a_i} \pm \sum_{j\neq
  k}\varepsilon'_{ij}a_jGf \pm \varepsilon'_{ik}a'_kGf)da_k =
$$
$$
=\int \pi i \hbar\widehat{d_i} G \frac{\partial f}{\partial a_i} \pm \sum_{j \neq
  k}\varepsilon_{ij}a_jGf \pm G \varepsilon_{ik}\left(a_k - \phi^{\hbar_k}\left(\sgn(\varepsilon_{ik})(\mp \pi i \hbar\widehat{d}_k\frac{\partial}{\partial a_k}-\sum_j\varepsilon_{kj}a_j)\right)f\right)da_k.
$$
Since these identities should be valid for any  $f$ one
gets the equations for the function $G$:
$$
\left(\pi i \hbar\widehat{d_i} \frac{\partial }{\partial
  a_i} \pm \sum_{j\neq k}(\varepsilon'_{ij}- \varepsilon_{ij})a_j\mp
\varepsilon_{ik}a_k \pm \varepsilon_{ik} \phi^{\hbar_k}\left(\sgn(\varepsilon_{ik})
(\pm \pi i \hbar\widehat{d}_k\frac{\partial}{\partial a_k}-
\sum_j\varepsilon_{kj}a_j)\right)\right)G=0.
$$
Let us introduce the function  $\widehat G$ related to $G$ by the 
Fourier transform:
$$
\widehat{G}(c)=\int e^{-\frac{a_kc}{\pi i \hbar}}G(a_k)da_k;~~
G(a_k)=\frac{1}{2\pi^2 \hbar}\int e^{\frac{a_kc}{\pi i \hbar}}\widehat{G}(c)dc
$$
(we omit the variables $a_1,\ldots,a_{k-1},a_{k+1},\ldots,\ldots,a_n$
both $G$ and $\widehat{G}$ depends on). 
Taking into account the relations
$$
\pi i \hbar \widehat{\frac{\partial G}{\partial
    a_k}}=c\widehat{G},~~\widehat{a_kG}=-\pi i \hbar\frac{\partial
  \widehat{G}}{\partial c}
$$
one can get the equation for the function $\widehat{G}$:
\begin{equation} \label{11.30.03.1}
\left(\pi i \hbar\widehat{d_i} \frac{\partial }{\partial
  a_i} \pm \sum_{j\neq k}(\varepsilon'_{ij}- \varepsilon_{ij})a_j\pm
\varepsilon_{ik}\pi i \hbar\frac{\partial}{\partial c} \pm 
\varepsilon_{ik} \phi^{\hbar_k}\left(\sgn(\varepsilon_{ik})(\pm 
\widehat{d}_kc-\sum_j\varepsilon_{kj}a_j)\right)\right)\widehat{G}=0.
\end{equation}
Taking sum and difference of the equations corresponding to the upper
and lower signs one obtains:
$$
\left(2 \pi i \hbar \widehat{d_i}\frac{\partial }{\partial a_i} + 
\varepsilon_{ik} \phi^{\hbar_k}\left(\sgn(\varepsilon_{ik})(\widehat{d}_kc -\sum_j\varepsilon_{kj}a_j)\right)- \varepsilon_{ik} \phi^{\hbar_k}
\left(\sgn(\varepsilon_{ik})(-\widehat{d}_kc -\sum_j\varepsilon_{kj}a_j)\right)\right)\widehat{G}=0
$$
and
$$
\left(2\sum_{j\neq k}(\varepsilon'_{ij}- \varepsilon_{ij})a_j+
2\varepsilon_{ik} \pi i h\frac{\partial}{\partial
  c}+\right.
$$
$$
\left.+\varepsilon_{ik} \phi^{\hbar_k}\left(\sgn(\varepsilon_{ik})(\widehat{d}_kc -\sum_j\varepsilon_{kj}a_j)\right)+
\varepsilon_{ik} \phi^{\hbar_k}\left(\sgn(\varepsilon_{ik})
(-\widehat{d}_kc -\sum_j\varepsilon_{kj}a_j)\right)\right)\widehat{G}=0.
$$
Observe that these are the system of $2n-2$ equations on a function of $n$
variables. So it is an overdetermined system if $n>2$. 
Using the identities
$$
\phi^{\hbar_k}(\sgn(a)b)=\phi^{\hbar_k}(b) + (\sgn(a)-1)b/2;~~\widehat{d}_i\varepsilon_{ji}=-\widehat{d}_j\varepsilon_{ij}
$$
they can be transformed to the form
\begin{equation} \label{12.1.03.2}
2\pi i \hbar \frac{\partial \log \widehat{G}}{\partial a_i} =\widehat{d}_k^{-1}\varepsilon_{ki} \left(\phi^{\hbar_k}(\widehat{d}_kc -\sum_j\varepsilon_{kj}a_j)-
\phi^{\hbar_k}(-\widehat{d}_kc -\sum_j\varepsilon_{kj}a_j) + \widehat{d}_kc(\sgn(\varepsilon_{ik})-1)\right)
\end{equation}
and
\begin{equation} \label{12.1.03.1}
2 \pi i \hbar \frac{\partial \log \widehat{G}}{\partial c}= 2(\varepsilon_{ik})^{-1}\sum_{j\neq
k}(\varepsilon_{ij}-\varepsilon'_{ij})a_j 
\end{equation}
$$
-\phi^{\hbar_k}(\widehat{d}_kc -\sum_j\varepsilon_{kj}a_j) - 
\phi^{\hbar_k}(- \widehat{d}_kc -\sum_j\varepsilon_{kj}a_j)+ 
(\sgn(\varepsilon_{ik})-1)\sum_j\varepsilon_{kj}a_j
$$
Taking into account that
$$
\varepsilon_{ij}'-\varepsilon_{ij}= \frac{
|\varepsilon_{ik}|\varepsilon_{kj}+\varepsilon_{ik}|\varepsilon_{kj}|}{2},
\qquad i,j \not = k
$$
we have the following identity:
$$
2(\varepsilon_{ik})^{-1}(\varepsilon_{ij}-\varepsilon'_{ij}) 
+ (\sgn(\varepsilon_{ik})-1)\varepsilon_{kj} = 
(\sgn(\varepsilon_{jk})-1)\varepsilon_{kj}
$$
Multiplying it by $a_j$ and taking the sum over $j 
\not =k$ we get
$$
2(\varepsilon_{ik})^{-1}\sum_{j\neq
  k}(\varepsilon_{ij}-\varepsilon'_{ij})a_j+(\sgn(\varepsilon_{ik})-1)\sum_j\varepsilon_{kj}a_j=\sum_j(\sgn(\varepsilon_{jk})-1)\varepsilon_{kj}a_j,
$$
Thus (\ref{12.1.03.1}) is equivalent to 
$$
2 \pi i \hbar \frac{\partial \log \widehat{G}}{\partial c}= 
-\phi^{\hbar_k}(\widehat{d}_kc -\sum_j\varepsilon_{kj}a_j) - 
\phi^{\hbar_k}(- \widehat{d}_kc -\sum_j\varepsilon_{kj}a_j) + 
(\sgn(\varepsilon_{jk})-1)\sum_j\varepsilon_{kj}a_j
$$
Therefore the solution of the equations (\ref{12.1.03.2}) and (\ref{12.1.03.1}) is given by the formula  
$$
\widehat{G}=C\Phi^{\hbar_k}(\widehat{d}_kc
-\sum_j\varepsilon_{kj}a_j)^{-1}
\Phi^{\hbar_k}(-\widehat{d}_kc
-\sum_j\varepsilon_{kj}a_j)e^{c\sum_j
(\sgn(\varepsilon_{jk})-1))\varepsilon_{kj}a_j/2
  \pi i \hbar},
$$
where $C$ is an arbitrary constant. Taking $C=2 \pi^2  \hbar$, one
obtains the desired formula (\ref{11.13.03.1}). 
The statement is proved. 

\begin{lemma}  \label{11.30.03.2}
An integral operator given by the formula (\ref{convolution}) for certain function $G$ 
  intertwines the operators 
(\ref{12.9.03.10}) if
and only if the standard formula  for the 
mutation of the function $\varepsilon_{ij}$ holds. 
\end{lemma}

{\bf Proof}. The proof of the theorem shows that this formula, 
 as well as the formula for mutations of the ${A}$-coordinates follow from 
the anzatz (\ref{convolution}) and the mutation formulas 
 for the quantized ${X}$-coordinates.

\vskip 3mm

{\it Representation of the modular double 
of the chiral double of ${\cal X}_{|{\bf i}|, q}$.}
 Combining Claim \ref{11.12.03.323} 
and Theorem \ref{11.11.03.76} 
we see that the collection of Hilbert spaces 
$\{L_2({\cal A}^+_{{\bf i}})\}$ should provide a projective unitary $\ast$-representation of the 
{\it modular double} 
of the chiral double of ${\cal X}_{|{\bf i}|, q}$, defined as 
$$
{\cal X}_{ |{\bf i}|, q} \times  
{\cal X}_{|{\bf i}^o|, q}
\times {\cal X}_{|{\bf i}^{\vee}|, {q^{\vee}}} \times 
{\cal X}_{ |{{\bf i}^o}^{\vee}|, {q^{\vee}}}.
$$

According to \cite{FG3},  the collection 
of Hilbert spaces 
$\{L_2({\cal A}^+_{{\bf i}})\}$ should provide a  representation of 
the modular double ${\cal D}_{ |{\bf i}|, {q}} \times 
{\cal D}_{|{\bf i}^{\vee}|, {q^{\vee}}}$ 
of the {\it cluster double} of the 
quantum cluster ${\cal X}$-variety  ${\cal X}_{ |{\bf i}|, {q}}$. 
Since there is a canonical map of quantum spaces
$
{\cal D}_{ |{\bf i}|, {q}} \lra {\cal X}_{ |{\bf i}|, {q}} \times  
{\cal X}_{|{\bf i}^o|, {q}}, 
$ 
this implies the above claim. 
\vskip 3mm

\section{The quantum logarithm and dilogarithm functions}

The proofs of all results listed above can be found in \cite{FG3}. 

\vskip 3mm
Recall the dilogarithm function
$$
{\rm Li}_2(x):= -\int_0^x\log(1-t)dt.
$$
{\bf The quantum logarithm function.} 
It is the following function:
\begin{equation} \label{phi}
\phi^\hbar(z):= - 2\pi \hbar\int_{\Omega}\frac{e^{-ipz}}{(e^{\pi p} - 
e^{-\pi p})
(e^{\pi \hbar p}-e^{-\pi \hbar p}) }dp; 
\end{equation}
where the contour $\Omega$ goes along the real axes 
from $- \infty$ to $\infty$ bypassing the origin 
from above. 

\begin{proposition} \label{5.31.03.1} The function $\phi^{\hbar}(x)$ 
enjoys the following properties. 
$$ \lim_{\hbar
\rightarrow 0}\phi^\hbar(z) = \log(e^z + 1).\leqno{\bf (A1)} 
$$
$$
\phi^\hbar(z)-\phi^\hbar(-z)=z. \leqno{\bf (A2)}  
$$
$$
\overline{\phi^\hbar(z)} = \phi^\hbar(\overline{z}).\leqno{\bf (A3)}
$$
$$
\phi^\hbar(z)/\hbar = \phi^{1/\hbar}(z/\hbar).\leqno{\bf (A4)}
$$
$$
\phi^\hbar(z+i\pi \hbar)-\phi^\hbar(z-i\pi \hbar) = \frac{2\pi i
\hbar}{e^{-z}+1}, \qquad 
 \phi^\hbar(z+i\pi)-\phi^\hbar(z-i\pi)
= \frac{2\pi i}{e^{-z/\hbar}+1}.\leqno{\bf (A5)}
$$
{\bf (A6)} The form
$\phi^\hbar(z)dz$ is meromorphic with poles at the points $\{\pi i ((2m-1)+
(2n-1)\hbar)|m,n \in {\mathbb N}\}$ with residues $2 \pi i \hbar$ and at the 
points $\{-\pi i ((2m-1)+(2n-1)\hbar)|m,n \in {\mathbb N}\}$ with residues 
$-2\pi i\hbar$.
\end{proposition}

{\bf The quantum dilogarithm.}
Recall the quantum dilogarithm function:
$$
\Phi^\hbar(z) := {\rm exp}\Bigl(-\frac{1}{4}\int_{\Omega}\frac{e^{-ipz}}{ {\rm sh} (\pi p)
{\rm sh} (\pi \hbar p) } \frac{dp}{p} \Bigr).
$$
It goes back to Barnes \cite{Ba}, and was used by Baxter \cite{Bax}, Faddeeev \cite{Fad}, 
and others. 

\begin{proposition} \label{qdilog} The function $\Phi^{\hbar}(x)$
enjoys the following properties.
$$
2 \pi i \hbar\,d \log \Phi^\hbar(z) =  \phi^\hbar(z)dz \leqno{\bf (B)}
$$
$$
\lim_{\Re z \rightarrow -\infty}\Phi^\hbar(z)=1.\leqno{\bf (B0)}
$$
Here the limit is taken along a line parallel to the real
axis. 
$$ \lim_{\hbar \rightarrow 0}\Phi^\hbar(z)/\exp \frac{-Li_2(-e^{z})}{2\pi i
  \hbar} = 1. \leqno{\bf (B1)}
$$
$$
\Phi^\hbar(z)\Phi^\hbar(-z)=\exp\left(\frac{z^2}{4\pi i \hbar}\right)
e^{-\frac{\pi i}{12}(\hbar+\hbar^{-1})}. \leqno{\bf (B2)}
$$
$$
\overline{\Phi^\hbar(z)} = (\Phi^\hbar(\overline{z}))^{-1}.\mbox{ In
  particular } |\Phi^\hbar(z)|=1 \mbox{ for } z\in {\mathbb R}.\leqno{\bf (B3)}
$$
$$
\Phi^\hbar(z) = \Phi^{1/\hbar}(z/\hbar).\leqno{\bf (B4)}
$$
$$
\Phi^\hbar(z+ 2 \pi i \hbar) = \Phi^\hbar(z) (1+qe^z), \qquad
\Phi^\hbar(z+ 2 \pi i ) = \Phi^\hbar(z) (1+ q^{\vee}e^{z/\hbar}).  \leqno{\bf (B5)}
$$
\noindent({\bf B6}) The function
$\Phi^\hbar(z)dz$ is meromorphic with poles at the points $$\{\pi i ((2m-1)+
(2n-1)\hbar)|m,n \in {\mathbb N}\}$$  and zeroes at the
points $$\{-\pi i ((2m-1)+(2n-1)\hbar)|m,n \in {\mathbb N}\}$$
\end{proposition}

\vskip 3mm
A.G.: Department of Math, Brown University, 
Providence RI 02906, USA; sasha@math.brown.edu

V.F.: ITEP, B. Cheremushkinskaya 25, 117259 Moscow, Russia.
fock@math.brown.edu

\end{document}